# Generalized replicator dynamics based on mean-field pairwise comparison dynamic

Hidekazu Yoshioka [a, *]

[a] Japan Advanced Institute of Science and Technology, 1-1 Asahidai, Nomi, Ishikawa 923-1292, Japan

[*] Corresponding author: yoshih@jaist.ac.jp, ORCID: 0000-0002-5293-3246

**Abstract**

The pairwise comparison dynamic is a forward ordinary differential equation in a Banach space whose solution is a time-dependent probability measure to maximize utility based on a nonlinear and nonlocal protocol. It contains a wide class of evolutionary game models, such as replicator dynamics and its generalization. We present an inverse control approach to obtain a replicator-type pairwise comparison dynamic from the large discount limit of a mean field game (MFG) as a coupled forward-backward system. This methodology provides a new interpretation of replicator-type dynamics as a myopic perception limit of the dynamic programming. The cost function in the MFG is explicitly obtained to derive the generalized replicator dynamics. We present a finite difference method to compute these models such that the conservation and nonnegativity of the probability density and bounds of the value function can be numerically satisfied. We conduct a computational convergence study of a large discount limit, focusing on potential games and an energy management problem under several conditions.

## 1. Introduction

### 1.1 Background

Social dynamics emerge because of dynamic decision-making among mutually interacting players. Evolutionary game models are pivotal mathematical tools for studying the social phenomena observed worldwide [1,2], such as water saving behavior [3], carbon emission trading [4], green production [5], traveler dynamics [6], urban road planning [7], research incentive design [8], and generic binary choice dynamics [9].

The governing equations of evolutionary game models are nonlinear ordinary differential equations (ODEs) in a Banach space [10]. The solutions to these equations are time-dependent probability measures that represent the optimized action profiles of the players, and their mathematical and numerical analyses are essential for assessing and predicting social dynamics. Most evolutionary game models are of the pairwise comparison type [11] where mutually interacting players determine their actions by comparing their utility values. This dynamic is known as pairwise comparison dynamic and contains various major existing evolutionary game models, such as replicator [12,13,14,15], logit [16,17,18], and projection [19,20] dynamics. The unique solvability and stability of pairwise comparison dynamics and related models have been previously studied [11,21].

A common drawback of evolutionary game models is the assumption that players determine their actions solely based on their current state. They are based on ODEs that evolve forward in time without incorporating any future predictions. In this sense, the players in evolutionary games are myopic. Mean field game (MFG) is another branch of dynamic game theory of mutually interacting players [22]. Players in the MFG make decisions such that an objective function, which is a conditional expectation containing future costs and utility, is optimized based on the current state, which is therefore less myopic than evolutionary game models in this sense. MFG is based on a dynamic programming principle that determines optimal actions from the future to the past using a time-backward evolution equation: the Hamilton–Jacobi–Bellman equation (HJB) [23,24]. The optimality equation of the MFG is a time forward-backward system of differential equations, with which both the forward evolution of the actions of players and their future predictions are concurrently simulated. Owing to its high flexibility, the MFG has been applied to various problems, including, but not limited to, microgrid energy management [25], smart energy charging [26,27], and mobility decisions during epidemics [28,29].

The evolutionary game model is derived based on pairwise microscopic interactions among players that recursively occur, and their aggregation appears as the governing equation [30]. Recently, studies have connected evolutionary games and MFGs to provide another explanation for the former and/or both. Bertucci et al. [31] discussed the convergence of an MFG system driven by diffusive players to a nonlinear Fokker–Planck (FP) equation. Bardi and Cardaliaguet [32] studied similar convergence problems in complex systems with degenerate diffusion. The central idea of their approach is that an evolutionary game model is obtained under an infinitely large discount rate in a proper MFG model such that the former is a myopic perception limit of the latter. From this perspective, the two different modeling frameworks are connected through a discount rate. Another stream of research exists in which an evolutionary game model is derived as the model predictive control version of some MFG [33-35]; they do not use the discount rate but discuss the convergence of the latter to the former by considering some myopic perception limit. Yoshioka and Tsujimura [36] and Yoshioka [37] argued that logistic-type evolutionary game models arise in the myopic perception limit of certain MFGs and presented computational studies concerning convergence speed. However, these studies do not cover the major classes of evolutionary game models, the replicator dynamics, which are frequently encountered in modern science and engineering problems [38-42]. Replicator-type dynamics are qualitatively different from logit ones in that the solution to the former does not expand its support [15], but the latter does. These backgrounds motivated the study explained in the next subsection.

### 1.2 Aim and contribution

The aim of this study is to present and investigate new linkages between a replicator-type dynamic, called the generalized replicator dynamic (GRD) in this paper, and the MFG, where the former turns out to be a myopic perception limit of the latter. With this new connection between the two seemingly different mathematical models obtained in this paper, it would become possible to understand and analyze dynamic games from a broader standpoint, advancing both theory and modelling technologies for dynamic decision-making in many directions.

The strategy taken in this paper is as follows:

Step 1: To formulate a GRD (**Section 2**)
Step 2: To formulate an MFG, and guess its relevant form based on the GRD (**Section 3**): direction from GRD to MFG
Step 3: To discuss that the GRD is indeed obtained as the large discount limit of the MFG (**Section 3**): direction from MFG to GRD
Step 4: To computationally validate the heuristics by choosing from small to large values of the discount rate (**Section 4**): numerically examine consistency between them

In this view, the MFG is designed based on GRD (Steps 1 and 2), and we demonstrate that it would be possible (Step 3) and computationally check our theory (Step 4). Our contributions to this research are explained in more detail as follows.

First, we formulate a GRD as a generalization of classical replicator dynamics such that its transition rate contains a wider class of nonlinearity. This formulation is performed within the framework of pairwise comparison dynamic. The well-posedness of the GRD is then established.

Second, we formulate an MFG in which the myopic perception limit is the GRD. Here, myopicity is represented by the discount of utility and cost in time [43,44]. The key technique is to formulate the control cost in the MFG based on an inverse control approach [45-48] to obtain a suitable objective function with a discount to be optimized in the game. The control cost, which is understood as the cost of revising the actions of players, is then linked to the transition rate of the GRD. A genetic formula for designing a cost function starting from a transition rate is derived, and specific cases with different nonlinearities, where the control cost is explicitly determined, are discussed. A heuristic discussion with proper scaling of coefficients in the MFG suggests that the myopic perception limit is the GRD. Thus, taking the limit of infinite impatience in an MFG leads to the GRD, an evolutionary game model. In this view, we can understand the GRD as a special limit of the MFG. Because an evolutionary game model is based on the infinite impatience, studying this limit is an important task.

Finally, we validate the theory for connecting the two dynamic game models through numerical computations based on a finite difference method. This finite difference method is structure-preserving in the sense that the bounds of numerical solutions, such as nonnegativity and some upper and lower bounds, are computationally preserved under certain conditions. Using this numerical method, we apply the GRD and its MFG counterpart to potential games and an energy management problem. In potential game cases, we consider both convex and concave utilities in a quadratic form. We demonstrate the convergence of the probability density and value function of the MFG to the probability density and utility in a GRD as the discount rate increases, suggesting that our formulation connects the two models through the discount rate. The energy management case assumes a more complex utility motivated by the sustainable use of hydropower to avoid degradation of aquatic environments. The convergence between the two formulations is again computationally observed. We also discuss the differences in the optimal actions of the two formulations under different transition rates and initial/terminal conditions.

Consequently, this study theoretically and computationally explores the linkage between evolutionary game and MFG, focusing on replicator-type dynamics. Our contributions provide new insights into the mathematical modeling of dynamic decision-making with a range of player foresight.

The remainder of this paper is organized as follows: **Section 2** reviews pairwise comparison dynamics and presents GRD, which is the main mathematical model studied in this paper. **Section 3** then derives the MFG related to the GRD based on an inverse control approach, showing that the two models are connected by the discount rate as a parameter. **Section 4** presents a computational analysis of the evolutionary and MFG models. **Section 5** concludes the study and presents perspectives. **Appendices** contain proofs of propositions and auxiliary results.

## 2. Generalized replicator dynamic

### 2.1 Notations

We prepare the notations to be used in this paper based on the literature on evolutionary games (e.g., Chapter 1.3 in Mendoza-Palacios and Hernández-Lerma [15]; Cheung et al. [11]). Time is a nonnegative parameter denoted as $t$. The domain of player actions is assumed to be a one-dimensional compact space and is given by $\Omega = [0,1]$ without loss of generality. The Borel sigma algebra of $\Omega$ is referred to as $\mathcal{B}$. The space of finite signed measures $\mu$ on $\mathcal{B}$ is denoted as $\mathcal{M}(\Omega)$ that is equipped with the total variation norm $\|\mu\| = \sup_g \left| \int_\Omega g(x)\mu(\mathrm{d}x) \right|$, where the supremum is taken for all measurable functions on $\Omega$ such that $|g| \leq 1$. Similarly, the space of probability measures on $\Omega$ is denoted as $\mathcal{P}(\Omega)$. A subspace of $\mathcal{M}(\Omega)$ with the total variation norm not greater than two ($\|\mu\| \leq 2$) is denoted as $\mathcal{M}_2(\Omega)$. An inclusion $\mathcal{P}(\Omega) \subset \mathcal{M}_2(\Omega) \subset \mathcal{M}(\Omega)$ exists. The time derivative of the time-dependent measure $\mu_t \in \mathcal{M}(\Omega)$ at each $t \geq 0$ is denoted as $\frac{\mathrm{d}\mu_t}{\mathrm{d}t}$, and is defined in a strong sense: $\lim_{\varepsilon \to 0} \left\| \frac{\mathrm{d}\mu_t}{\mathrm{d}t} - \frac{\mu_{t+\varepsilon} - \mu_t}{\varepsilon} \right\| = 0$. Unless otherwise specified, $A \in \mathcal{B}$ is an arbitrary Borel measurable set.

### 2.2 Pairwise comparison dynamic

A pairwise comparison dynamic is the following infinite-dimensional ODE in the form of a conservation law, where the temporal change in the probability of actions is determined by the balance between the inflow and outflow of probability currents [11]: given an initial condition $\mu_0 \in \mathcal{P}$, for $t > 0$,

$$\underbrace{\frac{\mathrm{d}}{\mathrm{d}t}\mu_t(A)}_{\text{Temporal change}} = \underbrace{\int_{x \in A}\int_{y \in \Omega} \rho(y, x, \mu_t)\mu_t(\mathrm{d}y)\mu_t(\mathrm{d}x)}_{\text{Inflow of probability}} - \underbrace{\int_{x \in A}\int_{y \in \Omega} \rho(x, y, \mu_t)\mu_t(\mathrm{d}y)\mu_t(\mathrm{d}x)}_{\text{Outflow of probability}}. \quad (1)$$

Here, $\rho : \Omega \times \Omega \times \mathcal{M}(\Omega) \to [0, +\infty)$ is the transition rate from some action to another. In (1), all players are assumed to be homogeneous and are characterized by actions in $\Omega$. A higher value of $\rho(x, y, \cdot)$ implies a larger probability flux of players from actions $x$ to $y$. Solution $\mu$ to dynamic (1) then describes temporal evolution of distribution of actions selected by players. Dynamic (1) is considered to be a nonlinear version of FP equations of a jump-driven stochastic process [49-51].

We specify the transition rate $\rho$ based on utility $U$ to be maximized by the players. Although the dynamic (1) is supposed to govern a probability measure whose total variation norm equals 1, utility should be defined in a wider space for technical reasons to define the dynamic well (e.g., Theorem 3.A in Zeidler [52]). In this paper, utility is a measurable function $U : \Omega \times \mathcal{M}(\Omega) \to \mathbb{R}$ that satisfies the following conditions: for any $x, y \in \Omega$ and $\mu, \nu \in \mathcal{M}_2(\Omega)$, a constant $K > 0$ exists with

$$|U(x, \mu)| \leq K \text{ (Boundedness)} \quad (2)$$

and

$$|U(x, \mu) - U(x, \nu)| \leq K\|\mu - \nu\| \text{ (Lipschitz continuity).} \quad (3)$$

A typical example complying with (2) and (3) is [10]

$$U(x, \mu) = \int_\Omega f(x, y)\mu(\mathrm{d}y) \quad (4)$$

with a bounded and measurable function $f : \Omega \times \Omega \to \mathbb{R}$.

### 2.3 Generalized replicator dynamic

The pairwise comparison dynamic (1) contains the classical replicator dynamic

$$\frac{\mathrm{d}}{\mathrm{d}t}\mu_t(A)=\int_{x\in A}\left(U(x,\mu_t)-\int_{y\in\Omega}U(y,\mu_t)\mu_t(\mathrm{d}y)\right)\mu_t(\mathrm{d}x). \tag{5}$$

Indeed, substituting the following $\rho$ into (1) yields (5):

$$\rho(x,y,\mu)=\left(U(y,\mu)-U(x,\mu)\right)_+, \tag{6}$$

where $(\cdot)_+=\max\{0,\cdot\}$. We generalize (5) and (6) where $\rho$ is given by

$$\rho(x,y,\mu)=C\left(U(y,\mu)-U(x,\mu)\right) \tag{7}$$

with some $C:\mathbb{R}\to[0,+\infty)$ that satisfies the following properties for any $x,y\in[0,+\infty)$: $C(0)=0$ (no transition when no utility difference is present), $C(x)>C(y)$ if $x>y$ (strictly increasing), $C(x)=0$ if $x<0$ (transit only when the utility after transition is larger than the current one), $C(x)<+\infty$ if $x<+\infty$ (boundedness), and $|C(x)-C(y)|\le L_{x,y}|x-y|$ with a constant $L_{x,y}>0$ depending on $x,y\in[0,+\infty)$ (local Lipschitz continuity). We investigate the following examples, where $q>0$ is a shape parameter and $x\ge 0$:

$$C(x)=x^q \text{ with } q\ge 1 \text{ (power type)}, \tag{8}$$

$$C(x)=\ln(qx+1) \text{ (logarithmic type)}, \tag{9}$$

$$C(x)=\exp(qx)-1 \text{ (positive exponential type)}, \tag{10}$$

$$C(x)=1-\exp(-qx) \text{ (negative exponential type)}. \tag{11}$$

The power type (8) reduces to the classical replicator dynamic (6) if $q=1$. Because a larger transition rate implies a more frequent transition for a large argument (i.e., a large utility difference), the transition is more likely to occur in the following order: exponential, power, logarithmic, and negative exponential types. Only the negative exponential type (11) admits a strictly bounded transition rate.

The GRD that we investigate in this study is the following evolution equation: given an initial condition $\mu_0\in\mathcal{P}(\Omega)$, for $t>0$,

$$\begin{aligned}\frac{\mathrm{d}}{\mathrm{d}t}\mu_t(A)=&\int_{x\in A}\int_{y\in\Omega}C\left(U(x,\mu_t)-U(y,\mu_t)\right)\mu_t(\mathrm{d}y)\mu_t(\mathrm{d}x)\\&-\int_{x\in A}\int_{y\in\Omega}C\left(U(y,\mu_t)-U(x,\mu_t)\right)\mu_t(\mathrm{d}y)\mu_t(\mathrm{d}x)\end{aligned}. \tag{12}$$

We conclude this subsection with the following proposition stating that dynamic (12) is well-posed.

***Proposition 1.***

*The GRD (12) admits a unique solution $(\mu_t)_{t\ge 0}$, that is, $\mu_t\in\mathcal{P}(\Omega)$ at each $t\ge 0$, and the left-hand side of (12) exists in the strong sense.*

### 2.4 Remarks

We present some remarks concerning GRD (12). Previous studies [36,37] investigated logit-type dynamics as unusual cases of pairwise comparison dynamics. A crucial difference between their evolutionary game models and the current studies is the support of probability measure. Here, the support of a probability measure is the maximal closed set $S \in \mathcal{B}$, such that $\mu(A) > 0$ for any $A \subset S$. A typical logit dynamic is as follows:

$$\frac{\mathrm{d}}{\mathrm{d}t}\mu_t(A) = \frac{\int_A \exp\left(U(x,\mu_t)/\eta\right)\mathrm{d}x}{\int_\Omega \exp\left(U(x,\mu_t)/\eta\right)\mathrm{d}x} - \mu_t(A), \tag{13}$$

where $\eta > 0$ represents the noise intensity during the decision-making process. The logit dynamic (13) forgets the support of the initial condition $\mu_0$ as time elapses and may formally approach the stationary state

$$\mu_t(A) = \frac{\int_A \exp\left(U(x,\mu_t)/\eta\right)\mathrm{d}x}{\int_\Omega \exp\left(U(x,\mu_t)/\eta\right)\mathrm{d}x} \tag{14}$$

whose right-hand side has the full support $S = \Omega$. By contrast, the pairwise comparison dynamic (1) preserves the support of initial condition $\mu_0$: if $\mu_0(A) = 0$ for some $A \in \mathcal{B}$, then $\mu_t(A) = 0$ at any $t > 0$. This implies that the stationary state of the dynamic, if it exists, depends on initial conditions. Moreover, it also implies that actions included in set $S$ will be never selected in the future. Indeed, this property was discussed for the replicator dynamic in Theorem 13 in Mendoza-Palacios and Hernández-Lerma [15]. The discussion here is heuristic, and a case in which stationary solutions to a logit dynamic depend on initial conditions exists [17]. Nevertheless, the difference between the logit and replicator dynamics would be particularly crucial if Nash equilibria are contained in the support of initial conditions or not. The support preservation property implies that GRD and its MFG counterpart presented in the next section should be investigated against different initial conditions having different supports. We investigate this issue in **Section 4**.

In the next section, based on the pairwise comparison dynamic reviewed above, we connect the GRD to an MFG by choosing a specific cost function. This connection between the two seemingly different mathematical models reveals that they can be understood from a unified viewpoint.

## 3. Mean field game

### 3.1 Formulation

#### 3.1.1 Stochastic model

We formulate an MFG counterpart of the GRD (12). The formulation in this section is based on a consistency condition and proceeds rather heuristically by assuming the existence of a large-population limit in an MFG because our focus is not on the rigorous microscopic derivation of an MFG but rather suggesting its linkage to the GRD [53-55].

We start from the single player case based on the stochastic calculus of jump processes under a complete probability space as usual (e.g., Øksendal and Sulem [56]). Consider a player whose action

$X=\left(X_t\right)_{t\ge0}$ follows the stochastic differential equation (SDE)

$$X_t=\begin{cases}\mathrm{d}L_t & (\text{at jump time of } L)\\ X_{t-} & (\text{otherwise})\end{cases},\ t>0, \tag{15}$$

where "$t-$" implies taking the left limit at time $t$, and $L=\left(L_t\right)_{t\ge0}$ is a compound Poisson process whose compensated measure at time $t$ is given by $u_{t-}(z)\mathrm{d}t\mathrm{d}z$ with $z$ being the jump size. Here, $u=\left(u_t\right)_{t\ge0}$ is a measurable mapping from $\Omega$ to $[0,+\infty)$ at each time $t$. The SDE (15) is equipped with an initial condition $X_0\in\Omega$. We assume $X_t\in\Omega$ almost surely for $t\ge0$. The jump intensity and distribution of this SDE are $\int_\Omega u_t(z)\mathrm{d}z$ and $u_t(z)/\int_\Omega u_t(z)\mathrm{d}z$ (if the denominator is nonzero), respectively. State $X$ is updated at jumps of $L$, and the size and distribution of jumps can be controlled by $u$.

We extend the single-player case to an $N$-player case with $N\in\mathbb{N}$, where action $X^{(i)}=\left(X_t^{(i)}\right)_{t\ge0}$ of the $i$ th player is governed by the SDE analogous to (15):

$$X_t^{(i)}=\begin{cases}\mathrm{d}L_t^{(i)} & (\text{at jump time of } L)\\ X_{t-}^{(i)} & (\text{otherwise})\end{cases},\ t>0 \tag{16}$$

subject to an initial condition $X_0^{(i)}\in\Omega$, where $L^{(i)}=\left(L_t^{(i)}\right)_{t\ge0}$ is a compound Poisson process with the compensated measure $t$ is $u_{t-}^{(i)}(z)\mathrm{d}t\mathrm{d}z$, where $u^{(i)}=\left(u_t^{(i)}\right)_{t\ge0}$ is a measurable mapping from $\Omega$ to $[0,+\infty)$ at each time $t$. We assume that jump times of each $L^{(i)}$ are based on mutually independent Poisson random measures.

Given $\left\{X_0^{(i)}\right\}_{i=1,2,3,\ldots,N}$ and $\left\{u^{(i)}\right\}_{i=1,2,3,\ldots,N}$, the objective function to be maximized by the $i$ th player is set as follows:

$$J_i\left(X_0^{(i)},u^{(i)}\right)=\mathbb{E}\left[\int_0^T e^{-\delta s}g\left(u_s^{(i)},X_s^{(i)},\mu_s^N\right)\mathrm{d}s+e^{-\delta T}\Psi\left(X_T^{(i)}\right)\right], \tag{17}$$

where $T>0$ is the prescribed terminal time, $\delta>0$ is the discount rate, $\mu_s^N$ is the empirical probability measure of $\left\{X^{(i)}\right\}_{i=1,2,3,\ldots,N}$ at time $s$, $\Psi:\Omega\to\mathbb{R}$ is the terminal gain, and $g:\Omega\times[0,+\infty)\times\mathcal{P}(\Omega)\to\mathbb{R}$ is the unit time cost and utility.

The MFG limit assumes a large-number limit $N\to+\infty$ under which the state dynamic is given by the SDE (15) whose control $u$ maximizes the objective function as follows:

$$\sup_u\phi\left(X_0,u\right)=\sup_u\mathbb{E}\left[\int_0^T e^{-\delta s}g\left(u_s,X_s,\mu_s\right)\mathrm{d}s+e^{-\delta T}\Psi\left(X_T\right)\right]. \tag{18}$$

The supremum is taken with respect to measurable and predictable (with respect to the natural filtration generated by $L$) time-dependent mappings from $\Omega$ to $[0,+\infty)$.

To apply the inverse control approach connecting the two different mathematical models, we design the coefficient $g$ as follows:

$$g(u_s, X_s, \mu_s) = \underbrace{\delta U(X_s, \mu_s)}_{\text{Scaled utility}} - \underbrace{\int_\Omega c(z, u_s)\,\mathrm{d}z}_{\text{Control cost}} \tag{19}$$

with some $c: \Omega \times [0, +\infty) \to [0, +\infty)$, which is fully specified later. In (19), the coefficient $g$ is separated into scaled utility and control costs, where the cost is assumed to arise when a player attempts to switch their actions. Scaling the utility term using the discount rate $\delta$ is crucial to our mathematical modeling of the myopic perception limit $\delta \to +\infty$.

For later use, the value function, as a parameterized version of (18) conditioned with information at time $t \in [0, T]$, is defined as follows:

$$\Phi(t, x) = \sup_u \mathbb{E}\left[\int_t^T e^{-\delta(s-t)}\left(\delta U(X_s, \mu_s) - \int_\Omega c(z, u_s)\,\mathrm{d}z\right)\mathrm{d}s + e^{-\delta(T-t)}\Psi(X_T)\middle| X_t = x\right]. \tag{20}$$

We assume that the optimal control $u = u^*$ that maximizes (20) exists and is Markovian.

### 3.1.2 Optimality equations

Assume that a density $\mu_t(\mathrm{d}x) = p_t(x)\mathrm{d}x$ exists. The dynamic programming principle (e.g., Cardaliaguet and Porretta [57]) suggests that the FP equation that governs $p$ is as follows:

$$\frac{\mathrm{d}}{\mathrm{d}t} p_t(x) = p_t(x)\int_{y \in \Omega}\left\{u_t^*(y, x) - u_t^*(x, y)\right\} p_t(y)\,\mathrm{d}y \tag{21}$$

for $0 < t \leq T$ and $x \in \Omega$ with a probability density as an initial condition, and the associated HJB equation to determine the value function $\Phi$ is as follows:

$$-\frac{\partial}{\partial t}\Phi(t, x) = -\delta\Phi(t, x) + \sup_u\left\{\int_\Omega \left(u(z)\left(\Phi(t, z) - \Phi(t, x)\right) - c(z, u(z))\right)\mathrm{d}z\right\} + \delta U(x, p_t) \tag{22}$$

for $0 \leq t < T$ and $x \in \Omega$ with the terminal condition $\Phi = \Psi$, where the supremum is taken with respect to the measurable function $u: \Omega \to [0, +\infty)$, and $u^*$ is a maximizer on the right-hand side of (22):

$$u_t^*(x, \cdot) = \underset{u(\cdot)}{\arg\max}\left\{\int_\Omega \left(u(z)\left(\Phi(t, z) - \Phi(t, x)\right) - c(z, u(z))\right)\mathrm{d}z\right\}. \tag{23}$$

Here, we express the utility through density $p$ instead of measure $\mu$ with an abuse of notation. According to the problem formulation, the argmax in (23) can be considered for the integrand of (22) at each $z \in \Omega$. The key to our modeling is to design the cost function $c$ such that the MFG system (21) and (22) effectively reduces to the GRD (12) under certain conditions.

## 3.2 Inverse control approach

The inverse control approach is a methodology for designing or estimating the cost function in a control problem such that the desired functional form of the optimal control can be analytically obtained [45-48]. This approach often requires a numerical algorithm to determine the cost function, whereas we analytically perform the task with the help of the assumed functional shape (7) of $\rho$. The dependence of the variables on $t$ is sometimes omitted for notational simplicity.

We set $\Delta_{zx} = \Phi(t,z) - \Phi(t,x)$. From (23), we expect that maximizing $u = u^*$ satisfies

$$\frac{\partial}{\partial u} c(z,u) = c'(z,u) = \Delta_{zx}, \tag{24}$$

where $c'(z,u)$ denotes the partial differential with respect to the second argument. If $\Phi(t,z) - \Phi(t,x) \geq 0$, we guess

$$u^* = c'^{(-1)}(z, \Delta_{zx}) = c'^{(-1)}(z, \Phi(t,z) - \Phi(t,x)), \tag{25}$$

where $c'^{(-1)}$ denotes the inverse of the second argument of $c'$. Considering the analogy between the FP equation (21) and the GRD (12), we expect the following:

$$c'^{(-1)}(z, \Phi(t,z) - \Phi(t,x)) = p(z) C(\Phi(t,z) - \Phi(t,x)), \tag{26}$$

where the value function $\Phi$ serves like a utility $U$ in the transition rate.

We analyze when (26) is satisfied. If $p(z) > 0$ at $z \in \Omega$, the answer is as follows: for any $u \geq 0$,

$$c(z,u) = \int_0^u C^{(-1)}\left(\frac{v}{p(z)}\right) \mathrm{d}v. \tag{27}$$

Indeed, if $u \geq 0$, then

$$\frac{\partial}{\partial u} c(z,u) = C^{(-1)}\left(\frac{u}{p(z)}\right). \tag{28}$$

If further $\Delta_{z,x} \geq 0$, using (24) and (28) we obtain

$$\Delta_{z,x} = C^{(-1)}\left(\frac{u^*}{p(z)}\right) \Leftrightarrow u^* = p(z) C(\Delta_{z,x}), \tag{29}$$

as desired. If $\Delta_{z,x} < 0$, then $u^* = 0$. Because of $C(x) = 0$ for $x \leq 0$, we arrive at the correct equality in (29) for any $\Delta_{z,x} \in \mathbb{R}$. If $p(z) = 0$, we can formally set $c(z,u) = +\infty$ to obtain $u^* = 0$.

Consequently, the desired cost function is derived as follows: if $u > 0$,

$$c(z,u) = \begin{cases} \int_0^u C^{(-1)}\left(\frac{v}{p(z)}\right) \mathrm{d}v & (p(z) > 0) \\ +\infty & (p(z) = 0) \end{cases} \tag{30}$$

and $c(z,u) = 0$ if $u = 0$. The designed cost $c(z,u)$ is increasing with respect to $u \geq 0$ and is decreasing with respect to $p(z)$. This cost therefore implies that jumping from the current action to the next one is more costly if the probability density of the next action is smaller, and vice versa. Imitating a major action, which is the one with a larger probability density, is less costly in this setting. Moreover, jumping to the action that no one is using is prohibitive. The designed cost of (30) thus inherits the imitative nature of replicator-type dynamics.

Considering (30), we have

$$\begin{aligned}\sup_u\left\{\int_\Omega\left(u(z)(\Phi(z)-\Phi(x))-c(z,u(z))\right)\mathrm{d}z\right\}&=\int_\Omega\left(p(z)C(\Delta_{z,x})\Delta_{z,x}-\int_0^{p(z)C(\Delta_{z,x})}C^{(-1)}\left(\frac{v}{p(z)}\right)\mathrm{d}v\right)\mathrm{d}z\\&=\int_\Omega p(z)\left(\int_0^{\Delta_{z,x}}C(w)\mathrm{d}w\right)\mathrm{d}z\end{aligned}\tag{31}$$

along with the following elementary identity when $p(z)>0$:

$$\begin{aligned}-\int_0^{p(z)C(\Delta_{z,x})}C^{(-1)}\left(\frac{v}{p(z)}\right)\mathrm{d}v&=-\int_0^{\Delta_{z,x}}C^{(-1)}\left(\frac{p(z)C(w)}{p(z)}\right)p(z)C'(w)\mathrm{d}w\\&=-p(z)\int_0^{\Delta_{z,x}}wC'(w)\mathrm{d}w\\&=-p(z)\left\{\left[wC(w)\right]_0^{\Delta_{z,x}}-\int_0^{\Delta_{z,x}}C(w)\mathrm{d}w\right\}\\&=-p(z)\Delta_{z,x}C(\Delta_{z,x})+p(z)\int_0^{\Delta_{z,x}}C(w)\mathrm{d}w\end{aligned}.\tag{32}$$

In summary, the HJB equation in the MFG becomes

$$-\frac{\partial}{\partial t}\Phi(t,x)=-\delta\Phi(t,x)+\int_\Omega p_t(z)\left(\int_0^{(\Phi(t,z)-\Phi(t,x))_+}C(w)\mathrm{d}w\right)\mathrm{d}z+\delta U(x,p_t).\tag{33}$$

The transition rate $C$ appears in the right-hand side of (33). Finally, the optimal transition rate $u^*$ under the MFG is determined as follows:

$$u_t^*=p_{t-}(\cdot)C\left(\Phi(t-,\cdot)-\Phi(t-,X_{t-})\right),\tag{34}$$

where the right-hand side provides a Markovian representation. The action after each jump is distributed according to the probability density just before the jump as well as the difference in value functions; the latter encodes some future predictions according to a conditional expectation (20).

We rewrite (33) as follows:

$$\Phi(t,x)-U(x,p_t)=\frac{1}{\delta}\left\{\frac{\partial}{\partial t}\Phi(t,x)+\int_\Omega p_t(z)\left(\int_0^{(\Phi(t,z)-\Phi(t,x))_+}C(w)\mathrm{d}w\right)\mathrm{d}z\right\}.\tag{35}$$

If $\Phi$ is bounded irrespective of $\delta$, which is justified at least computationally as shown in **Proposition 2**, letting $\delta\to+\infty$ during the intermediate time interval would yield $\Phi(t,x)=U(x,p_t)$ because the right-hand side of (35) vanishes. This implies that, at $0<<t<<T$ under $\delta\to+\infty$, we have $\Phi(t,x)=U(x,p_t)$, as along with the proximity of the probability densities between the GRD and MFG. In this view, the MFG heuristically leads to the GRD, providing a new derivation of replicator-type dynamics. Consequently, we demonstrated that the GRD can be derived as a large discount limit of the MFG, suggesting that the two models can be understood in a unified viewpoint by seeing the discount rate $\delta$ as a parameter to characterize these models. A similar reasoning has been successfully employed for logit dynamics [36,37]. We computationally demonstrate that this ansatz holds in **Section 4**.

Finally, for cases (8)–(11), the inner integral on the right-hand side of (33) can be explicitly determined as follows: These formulae are useful for numerically computing the MFG where we assume

$\Delta_{z,x} \geq 0$:

$$\int_0^{\Delta_{z,x}} C(w)\mathrm{d}w = \frac{1}{q+1}\Delta_{z,x}^{q+1} \text{ (power type)}, \tag{36}$$

$$\int_0^{\Delta_{z,x}} C(w)\mathrm{d}w = \frac{(q\Delta_{z,x}+1)\ln(q\Delta_{z,x}+1) - q\Delta_{z,x}}{q} \text{ (logarithmic type)}, \tag{37}$$

$$\int_0^{\Delta_{z,x}} C(w)\mathrm{d}w = \frac{\exp(q\Delta_{z,x}) - q\Delta_{z,x} - 1}{q} \text{ (positive exponential type)}, \tag{38}$$

$$\int_0^{\Delta_{z,x}} C(w)\mathrm{d}w = \frac{\exp(-q\Delta_{z,x}) + q\Delta_{z,x} - 1}{q} \text{ (negative exponential type)}. \tag{39}$$

### 3.3 Discussion

The derivation procedure for the cost function in the MFG assumes that a density $p$ for the probability measure $\mu$ exists. Nevertheless, one formally has

$$\int_\Omega p_t(z)\left(\int_0^{(\Phi(t,z)-\Phi(t,x))_+} C(w)\mathrm{d}w\right)\mathrm{d}z = \int_\Omega \left(\int_0^{(\Phi(t,z)-\Phi(t,x))_+} C(w)\mathrm{d}w\right)\mu_t(\mathrm{d}z) \tag{40}$$

whose right-hand side is well defined if $\Phi(t,z)-\Phi(t,x)$ is bounded and measurable. The derivation procedure, particularly the cost function representation in (30), should then be justified for $\mu$.

Given a time-dependent probability density $p = (p_t)_{t\geq 0}$, the HJB equation can be rewritten using a Hamiltonian $H : \Omega \times \mathcal{P}(\Omega) \times \mathbb{R} \times B(\Omega) \to \mathbb{R}$ as ($B(\Omega)$ is a space of bounded functions from $\Omega$ to $\mathbb{R}$) as follows:

$$-\frac{\partial}{\partial t}\Phi(t,x) = H\left(x, p_t, \Phi(t,x), \Phi(t,\cdot) - \Phi(t,x)\right) \tag{41}$$

with

$$H(x, p_t, \alpha, \beta) = -\delta\alpha + \int_\Omega p_t(z)\left(\int_0^{(\beta)_+} C(w)\mathrm{d}w\right)\mathrm{d}z + \delta U(x, p_t). \tag{42}$$

Because $C$ is nonnegative and strictly increasing, the quantity $\int_0^{(\Phi(t,z)-\Phi(t,x))_+} C(w)\mathrm{d}w$ is nonnegative, nondecreasing, and locally Lipschitz continuous with respect to $\Phi(t,z)-\Phi(t,x)$; for any $\Delta \in \mathbb{R}$, we have

$$\frac{\mathrm{d}}{\mathrm{d}\Delta}\int_0^{(\Delta)_+} C(w)\mathrm{d}w = C(\Delta) \times \begin{cases} 1 & (\Delta \geq 0) \\ 0 & (\Delta < 0) \end{cases}. \tag{43}$$

This implies that the Hamiltonian $H$ is nonincreasing and nondecreasing with respect to the second and third arguments, respectively. This is a condition known as nonlocal degenerate ellipticity, which is essential for defining HJB-type equations in viscosity sense [57-59]. If utility $U$ is continuous with respect to the first argument and the primitive function of $C$ is globally Lipschitz continuous, then the HJB equation (41) which is considered to be a single nonlocal partial differential equation, admits at most one viscosity solution that is continuous both in space and time (e.g., Theorem 3 Barles and Imbert [60]) with a modification that the present HJB equation is defined both inside a domain and its boundary). This solution can be obtained as a converged sequence of numerical solutions to the corresponding finite difference method, which is monotone, stable, and consistent (see **Propositions 2** presented later) [61,61].

Our MFG is of the non-separable type, where the probability density $p$ and value function $\Phi$ jointly appear in a single term, which is the integral term in (33). Problems of a non-separable nature are more difficult to address because limited techniques are applied to their numerical computation [63]. The non-separable nature poses another difficulty, such that the strict monotonicity is not satisfied or is not straightforward to check. A requirement for the monotonicity in our context is expressed as (e.g., Assumption 2 in Gomes et al. [64]): a constant $\gamma > 0$ exists such that

$$\frac{\partial^2}{\partial \Delta^2} \int_0^{(\Delta)_+} C(w) \mathrm{d}w \geq \gamma \quad \text{for all } \Delta \in \mathbb{R} . \tag{44}$$

However, this inequality is not satisfied because $\int_0^{(\Delta)_+} C(w) \mathrm{d}w = 0$ for $\Delta < 0$ with which the left-hand side of (44) equals zero.

Finally, we investigate a continuous-time setting, but its formulation is carried over to a discrete-time setting with suitable modifications. The discretized governing equations presented in the next section can be seen a discrete-time model [65].

## 4. Computational investigations

To answer the question of how to design policies considering his/her foresights (the discount rate in this paper) would be a research goal of our study. For example, the planning of resource or energy would be based on some future predictions dependent on policy-makers. Although partially, we try to answer this question in this section computationally.

We present and investigate a finite difference method (Section 4.1) to compute the MFG and GRD. With this numerical method, we analyze the two models computationally. More specifically, we have two goals in this section. We demonstrate that the MFG converges to GRD as the discount rate increases in potential game cases as the simplest problems in game theory (Section 4.2), the first goal. This computationally supports the heuristics discussed in Section 3. We also examine the convergence of the MFG to GRD in a more complex problem motivation from an application (Section 4.3), the second goal. Moreover, we quantify the speed of convergence from MFG to GRD to better understand their connection, which was absent in the theoretical analysis.

### 4.1 Finite difference method

We numerically discretize the forward-backward system (21) and (33). The space-time computational grid is the collection of points $\mathrm{P}_{i,j} : (t_i, x_j) \in [0,T] \times \Omega$ with $t_i = i\Delta t$ ( $i = 0,1,2,...,I$ ) and $x_j = (j-1/2)\Delta x$ ( $j = 1,2,3,...,J$ ), where $\Delta t = T/I$, $\Delta x = 1/J$, and $I, J \in \mathbb{N}$. The discretized quantity at $(t_i, x_j)$ is represented by subscript " $i, j$ ." Our numerical method is a finite difference method [37,66,67]. We apply an explicit Euler discretization to both the FP and HJB equations such that they are temporally integrated without resorting to the inversion of any matrix systems. The discretization of the FP equation is conservative, whereas that of the HJB equation preserves the upper and lower bounds of the solutions as shown later.

The discretization at $\mathrm{P}_{i,j}$ is as follows: FP equation

$$\frac{p_{i,j} - p_{i-1,j}}{\Delta t} = p_{i-1,j}\Delta x \sum_{k=1}^{J} \left( u_{i,k,j}^{*} - u_{i,j,k}^{*} \right) p_{i-1,k} \tag{45}$$

for $i = 1,2,3,...,I$ and $j = 1,2,3,...,J$ subject to an initial condition $\{p_{0,j}\}_{j=1,2,3,...,J}$ : HJB equation

$$\frac{\Phi_{i,j} - \Phi_{i+1,j}}{\Delta t} = -\delta \Phi_{i+1,j} + \sum_{k=1}^{J} p_{i,k} \left( \int_{0}^{(\Phi_{i+1,k} - \Phi_{i+1,j})_{+}} C(w)\,\mathrm{d}w \right) \Delta x + \delta U_{i,j} \tag{46}$$

for $i = 0,1,2,...,I-1$ and $j = 1,2,3,...,J$ , subject to the terminal condition $\Phi_{I,j} = \Psi(x_j)$ for $j = 1,2,3,...,J$ : optimal transition rate

$$u_{i,j,k}^{*} = C\left(\Phi_{i,k} - \Phi_{i,j}\right), \tag{47}$$

for $i = 0,1,2,...,I$ and $j,k = 1,2,3,...,J$ : utility

$$U_{i,j} = U\left(x_j, \{p_{i,k}\}_{k=1,2,3,...,J}\right) \tag{48}$$

for $i = 0,1,2,...,I$ and $j = 1,2,3,...,J$ . For the utility of form (4), we use

$$U_{i,j} = \Delta x \sum_{k=1}^{J} f\left(x_j, x_k\right) p_{i,k} . \tag{49}$$

The discretization for the GRD is based on (45):

$$\frac{p_{i,j} - p_{i-1,j}}{\Delta t} = p_{i-1,j}\Delta x \sum_{k=1}^{J} \left( C\left(U_{i-1,j} - U_{i-1,k}\right) - C\left(U_{i-1,k} - U_{i-1,j}\right) \right) p_{i-1,k} \tag{50}$$

for $i = 1,2,3,...,I$ and $j = 1,2,3,...,J$ subject to an initial condition $\{p_{0,j}\}_{j=1,2,3,...,J}$ .

We have two propositions, implying that the proposed numerical method generates reasonable numerical solutions that comply with the properties of the original problems.

***Proposition 2 (Stability of discretized value function)***

*Assume that the primitive function of* $C$ *is globally Lipschitz continuous with the Lipschitz constant* $\bar{L} > 0$ . *Moreover, assume* $(\delta + L_C)\Delta t < 1$ . *Given a nonnegative sequence* $\{p_{i,j}\}_{\substack{i=0,1,2,...,I \\ j=1,2,3,...,J}}$ *such that* $\sum_{j=1}^{J} p_{i,j}\Delta x = 1$ ( $i = 0,1,2,...,I$ ), *any sequence* $\{\Phi_{i,j}\}_{\substack{i=0,1,2,...,I \\ j=1,2,3,...,J}}$ *that solves (46) satisfies* $|\Phi_{i,j}| \le K_2$

*(* $i = 0,1,2,...,I$ *and* $j = 1,2,3,...,J$ *) with a constant* $K_2 > 0$ *independent of* $\Delta t$*,* $\Delta x$*, and* $\delta$*.*

***Proposition 3 (Conservation, nonnegativity, support preservation of probability density)***

*Assume that* $2C(2K_2)\Delta t < 1$*. Given a nonnegative sequence* $\{\Phi_{i,j}\}_{\substack{i=0,1,2,...,I \\ j=1,2,3,...,J}}$*, such that* $|\Phi_{i,j}| \le K_2$ *(* $i = 0,1,2,...,I$ *and* $j = 1,2,3,...,J$ *) with a constant* $K_2 > 0$ *independent of* $\Delta t$*,* $\Delta x$*. Then, any sequence* $\{p_{i,j}\}_{\substack{i=0,1,2,...,I \\ j=1,2,3,...,J}}$ *that solves (45) satisfies the following properties:*

*(2.a)* $$\sum_{j=1}^{J} p_{i,j}\Delta x = 1 \ (i = 0,1,2,...,I), \quad (51)$$

*(2.b)* $$p_{i,j} \ge 0 \ \text{if} \ p_{0,j} \ge 0 \ (i = 0,1,2,...,I, \ j = 1,2,3,...,J), \quad (52)$$

*(2.c)* $$p_{i,j} = 0 \ \text{if} \ p_{0,j} = 0 \ (i = 0,1,2,...,I) \ \text{for some} \ j = 1,2,3,...,J. \quad (53)$$

System (50) simply evolves in a time-forward manner, whereas the discretized MFG system (45) and (46) is a forward-backward system that should be handled using an iteration method. We use a fixed-point iteration method that has been applied to several engineering problems [68-70]. We used a relaxation factor of 0.25 such that the convergence of the numerical solutions is achieved in a more stable manner. The number of fixed-point iterations in our setting in **Section 4** was approximately 100 to 200.

The assumption of global Lipschitz continuity of $C$ seems severe; however, we can truncate $C$ using a sufficiently large constant when its argument is large, for instance, $3K_2$ because $K_2$ does not depend explicitly on $C$ (see **Proof of Proposition 3**). Even with this truncation, the power case with $m \in (0,1)$ should be excluded because of its non-Lipschitz continuity at the origin. Consequently, we have the following corollary.

***Corollary 1 (Removal of the global Lipschitz continuity assumed in Proposition 3)***

*The conclusion of* ***Proposition 2*** *holds true when the primitive function of* $C$ *is Lipschitz continuous in the compact interval* $[0, 3K_2]$*.*

By **Propositions 2 and 3**, the convergence of the numerical solution to the FP and HJB equations, considered to be a single equation, can be established in measure (Chapter 5 in Mendoza-Palacios and Hernández-Lerma [15]) and viscosity senses (Yoshioka and Tsujimura [36]), respectively, while the issue of the coupled forward-backward system would be more complicated. Finally, implicit schemes analogous to existing ones [71,72] can also be implemented if preferred. In this study, we simplify the numerical method and will examine these schemes in the future.

### 4.2 Potential games

The first application contains the potential cases discussed in Friedman and Ostrov [73] where the utility is given by (4) with either $f(x,y) = -(x-y)^2$ (concave case) or $f(x,y) = (x-y)^2$ (convex case). A Nash equilibrium is the Dirac delta concentrated at 1/2 (concave case) and the sum of Dirac deltas concentrated at 0 and 1 with equal weights 1/2 (convex case). Concave $f$ corresponds to a potential game with a concave potential, and hence it admits a monotone utility. By contrast, convex $f$ corresponds to a potential game with a convex potential, and hence it admits an anti-monotone utility. This difference suggests larger difficulty in stably obtaining numerical solutions in the latter case. We specify the computational resolution $I = 10000$ and $J = 200$. This computational resolution is satisfactory for our application as shown in **Section A.2 of the Appendix**. Moreover, the computational error is proportional to the reciprocal of $I, J$, suggesting a first-order convergence in the finite difference method.

First, we examine the uniform initial distribution of the probability density $p$. The terminal time for the MFG is $T = 100$, and the terminal condition is $\Psi = 0$. We compare the behaviors of the probability density, utility, and value functions for different coefficients $C$ and different values of the discount rate $\delta$. **Table 1** (power case with $q = 1$, classical replicator case) **and Table 2** (exponential case with $q = 2$) list the maximum absolute errors between the probability densities of the GRD and MFG for concave $f$. The comparison results of the utility of the GRD and the value function of the MFG are also summarized in **Tables 1 and 2**. Convergence from the MFG to GRD is achieved at the speed of $O(\delta^{-1})$, supporting our heuristic argument that the latter is a myopic perception limit in terms of both probability density and utility. Similarly, **Table 3** (power case with $q = 1$, classical replicator case) and **Table 4** (exponential case with $q = 2$) summarize the results for convex $f$. In this case, convergence from the MFG to the GRD is again achieved at the rate of $O(\delta^{-1})$; however, convergence failures occur for a small $\delta$. The multi-peak nature of the probability density is considered to degrade the convergence of the numerical algorithm, or the well-pawedness of the problem itself has been broken owing to convex utility, as discussed in **Section 3.3**. The different coefficients $C$ in (8)–(11) qualitatively achieve the same convergence rate (CV), except for a small $\delta$ for convex $f$ where the convergence fails (see **Section A2 in the Appendix**).

The convergence discussed above is visualized in **Figure 1** for the probability density and **Figure 2** for utility (value function) for the classical replicator case, again showing the convergence of the MFG to the GRD as $\delta$ increases, except near the terminal time $T$, supporting the convergence. **Figure 3** for the exponential case ($q = 2$) shows the corresponding results along the centerline $x = 0.5$. The computational results for the logarithmic and negative exponential cases are presented in **Section A.3 of the Appendix** and further support the consistency between the two models.

We also analyze a finite-supported initial condition with $p(x) = 5/3$ when $x \le 3/5$, and $p(x) = 0$ otherwise. **Tables 5 and 6** suggest the convergence of the MFG to the GRD even under this restricted initial condition, implying its wide applicability.

**Table 1.** Comparison of the probability density and value function (utility in GRD) between the classical replicator dynamic and corresponding MFG for concave $f$. CV at $\delta = 10^{m}$ against $\delta = 10^{m-1}$ is evaluated as $\log_{10}\left(\mathrm{Err}_{\delta=10^{m-1}} / \mathrm{Err}_{\delta=10^{m}}\right)$, where Err represents the absolute error between the two probability densities. The Err in the column "Value function" is evaluated by the difference between the utility of GRD and value function of the MFG. The same applies to the other tables presented later.

| | Probability density | | Value function | |
|---|---|---|---|---|
| $\delta$ | Err | CV | Err | CV |
| $10^{-2}$ | 1.4.E+00 | | 1.9.E-01 | |
| $10^{-1}$ | 1.5.E-01 | 9.64.E-01 | 9.8.E-02 | 2.79.E-01 |
| $10^{0}$ | 1.5.E-02 | 9.93.E-01 | 2.3.E-02 | 6.29.E-01 |
| $10^{1}$ | 1.5.E-03 | 9.99.E-01 | 2.8.E-03 | 9.22.E-01 |
| $10^{2}$ | 1.6.E-04 | 9.78.E-01 | 2.8.E-04 | 9.93.E-01 |

**Table 2.** Same as **Table 1** for the exponential case ($q = 2$).

| | Probability density | | Value function | |
|---|---|---|---|---|
| $\delta$ | Err | CV | Err | CV |
| $10^{-2}$ | 2.0.E+00 | | 2.0.E-01 | |
| $10^{-1}$ | 2.1.E-01 | 9.78.E-01 | 1.3.E-01 | 1.82.E-01 |
| $10^{0}$ | 2.2.E-02 | 9.94.E-01 | 4.5.E-02 | 4.69.E-01 |
| $10^{1}$ | 2.1.E-03 | 1.00.E+00 | 6.6.E-03 | 8.35.E-01 |
| $10^{2}$ | 2.2.E-04 | 9.90.E-01 | 6.9.E-04 | 9.77.E-01 |

**Table 3.** Same as **Table 1** for convex $f$. "CF" denotes convergence failure; the same applies to tables presented later.

| | Probability density | | Value function | |
|---|---|---|---|---|
| $\delta$ | Err | CV | Err | CV |
| $10^{-2}$ | CF | | CF | |
| $10^{-1}$ | CF | | CF | |
| $10^{0}$ | 1.4.E-01 | | 2.2.E-02 | |
| $10^{1}$ | 1.4.E-02 | 1.00.E+00 | 2.6.E-03 | 9.25.E-01 |
| $10^{2}$ | 1.4.E-03 | 1.01.E+00 | 2.6.E-04 | 9.97.E-01 |

**Table 4.** Same as **Table 2** for convex $f$.

| | Probability density | | Value function | |
|---|---|---|---|---|
| $\delta$ | Err | CV | Err | CV |
| $10^{-2}$ | CF | | CF | |
| $10^{-1}$ | CF | | CF | |
| $10^{0}$ | 2.5.E-01 | | 4.4.E-02 | |
| $10^{1}$ | 2.5.E-02 | 1.00.E+00 | 6.4.E-03 | 8.36.E-01 |
| $10^{2}$ | 2.5.E-03 | 9.98.E-01 | 6.8.E-04 | 9.78.E-01 |

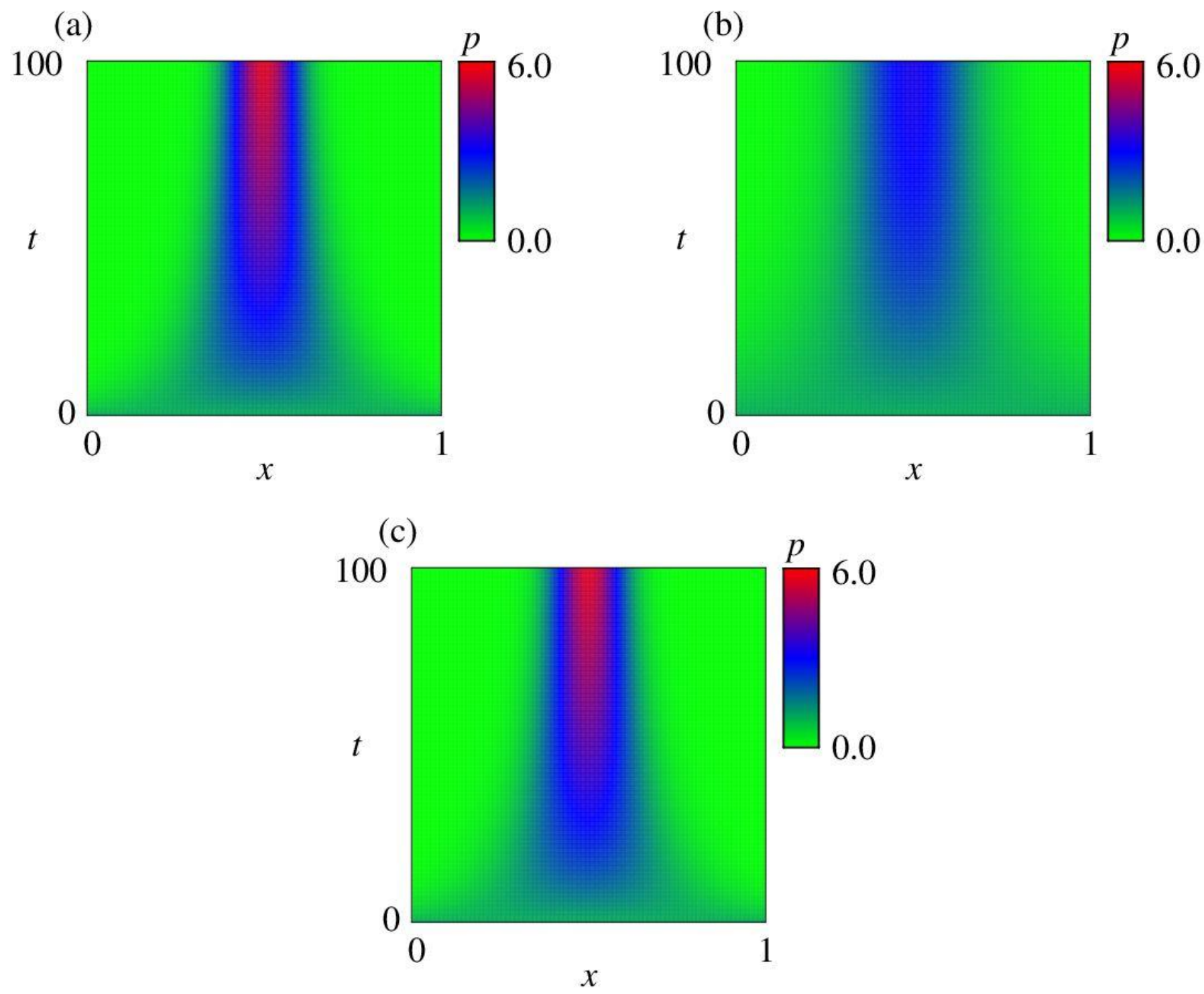

**Figure 1.** Comparison of the probability densities $p$ of the classical replicator case: the (a) GRD, (b) MFG with $\delta = 0.01$, and (c) MFG with $\delta = 1$.

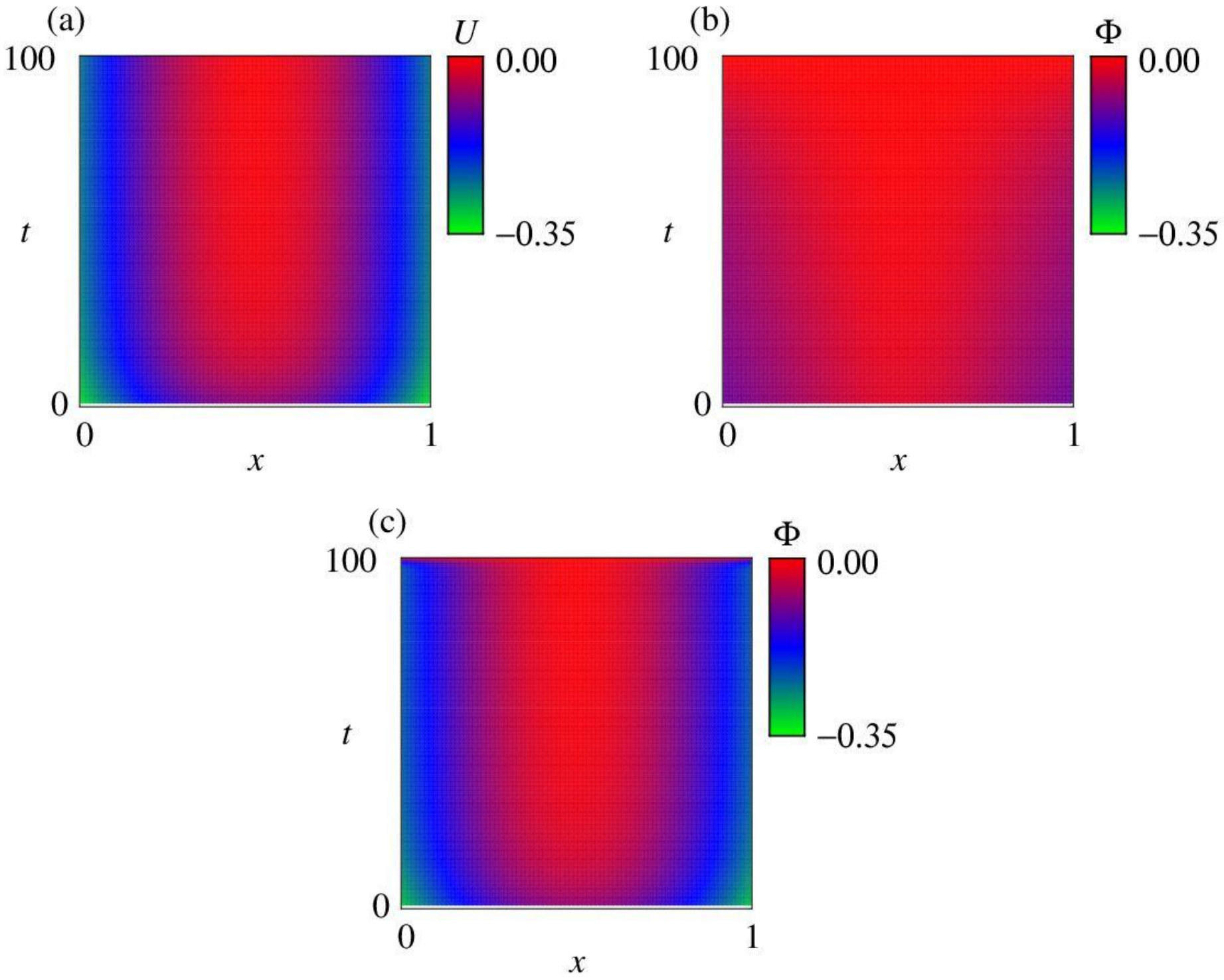

**Figure 2.** Comparison of utility $U$ and value function $\Phi$ for the classical replicator case: (a) utility of the GRD, (b) value function of the MFG with $\delta = 0.01$, and (c) value function of the MFG with $\delta = 1$.

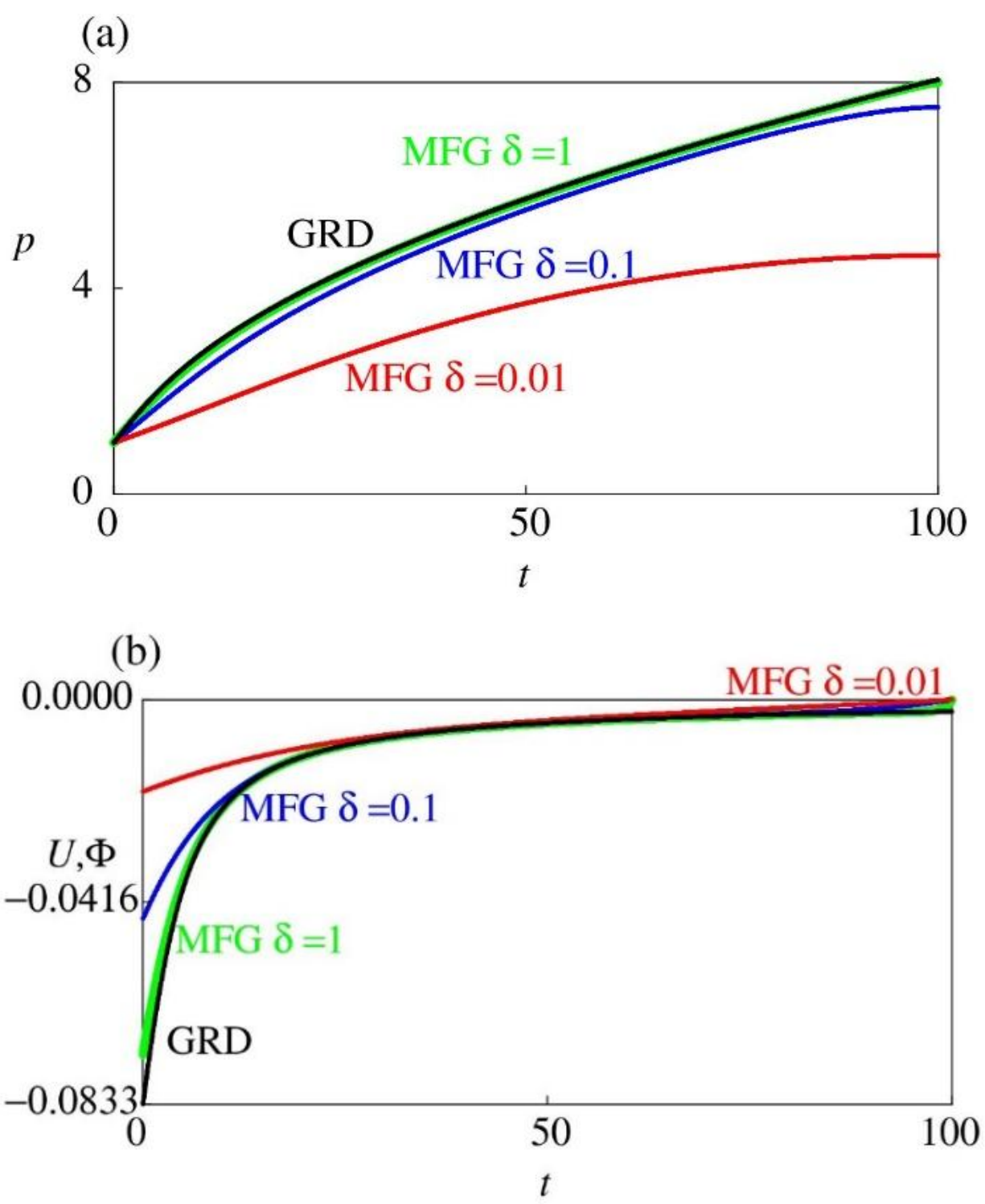

**Figure 3.** Comparison of (a) the probability density $p$ and (b) utility $U$ and value function $\Phi$ between the GRD and MFG for different values of $\delta$ for the exponential case ( $q = 2$ ).

**Table 5.** Same as **Table 1** for the exponential case ( $q = 2$ ) and concave $f$ but with the finite-support initial condition.

| | Probability density | | Value function | |
|---|---|---|---|---|
| $\delta$ | Err | CV | Err | CV |
| $10^{-2}$ | 1.9.E+00 | | 4.2.E-01 | |
| $10^{-1}$ | 2.1.E-01 | 9.54.E-01 | 3.1.E-01 | 1.27.E-01 |
| $10^{0}$ | 2.2.E-02 | 9.89.E-01 | 1.4.E-01 | 3.38.E-01 |
| $10^{1}$ | 2.2.E-03 | 1.00.E+00 | 2.8.E-02 | 7.05.E-01 |
| $10^{2}$ | 2.2.E-04 | 1.00.E+00 | 3.2.E-03 | 9.45.E-01 |

**Table 6.** Same as **Table 1** for the logarithmic case ( $q = 0.5$ ) and concave $f$ but with the finite-support initial condition.

| | Probability density | | Value function | |
|---|---|---|---|---|
| $\delta$ | Err | CV | Err | CV |
| $10^{-2}$ | 6.3.E-01 | | 3.7.E-01 | |
| $10^{-1}$ | 7.9.E-02 | 8.96.E-01 | 1.9.E-01 | 2.93.E-01 |
| $10^{0}$ | 8.9.E-03 | 9.50.E-01 | 4.3.E-02 | 6.48.E-01 |
| $10^{1}$ | 9.0.E-04 | 9.95.E-01 | 5.1.E-03 | 9.30.E-01 |
| $10^{2}$ | 9.0.E-05 | 1.00.E+00 | 5.1.E-04 | 9.94.E-01 |

### 4.3 Energy management application

We consider a simple common-resource management problem for energy users, whose energy use is represented by a normalized variable in $\Omega$ . The utility in this problem is given by

$$U(x,\mu) = \underbrace{\frac{1}{\alpha}x^{\alpha}}_{\text{Utility by energy use}} - \underbrace{\sigma\left(1 + w\underbrace{\int_{[\bar{X},1]}\mu(\mathrm{d}y)}_{\text{Density effect}}\right)x}_{\text{Energy cost}} \tag{54}$$

where $\alpha \in (0,1)$ (shape parameter of utility), $\sigma > 1$ (unit cost of energy use), $w > 0$ (proportional coefficient of the density effect), and $\bar{X} \in (0,1]$ (threshold above which the density effect is activated). This utility models the concave utility of energy use subject to a nonlinear cost of energy use, where the cost can be designed by the energy supplier considering the distribution of energy use by users such that the energy source as a common-pool resource is not extracted with an intensity higher than $\bar{X}$ .

The present case is fundamentally different from the potential ones in the previous subsection because it is not a potential game owing to asymmetry; indeed, utility (54) is rewritten in the form (4) with $f$ that satisfies $f(x,y) \neq f(y,x)$:

$$f(x,y) = \frac{1}{\alpha}x^{\alpha} - \sigma x - \sigma w \begin{cases} 0 & (0 \le x < \bar{X}) \\ 1 & (\bar{X} \le x \le 1) \end{cases}. \tag{55}$$

This utility complies with conditions (2) and (3); hence, the GRD is well-posed.

Utility (54) motivates hydropower generation, where the excessive intake of river water for electricity supply leads to the degradation of the water environment owing to extremely low water levels. A threshold water level below which the life of aquatic organisms will be threatened exists [74,75]. The energy cost can then be designed such that it will be increased to suppress energy use by users if it affects the water environment, that is if $\int_{[\bar{X},1]} \mu(\mathrm{d}y) > 0$, where hydropower is generated. As a more detailed example, the Tedori River Dam in Ishikawa Prefecture, Japan (dam data), has a hydropower generation system that supplies electricity to $O(10^5)$ residents and regulates the outflow, varying from 0 ($m^3$/s) to 200 ($m^3$/s) with an average of 34 ($m^3$/s) and a standard deviation of 40 ($m^3$/s) (the data from April 1 2016 to March 31 2024)[1]. The low flow, implying almost no flow in the downstream river, is attributed to the intake of water for hydropower generation. This type of flow regulation would critically affect local aquatic organisms and can be effectively mitigated if water use is decreased through policymaking by suppliers. The second term of utility (54) is intended to avoid this situation.

The expected equilibria are briefly reported and are investigated in detail in **Section A.3 of the Appendix** because the focus here is the convergence between the two models. Set $\bar{X}_1 = \sigma^{-\frac{1}{1-\alpha}}$ and $\bar{X}_2 = (\sigma(1+w))^{-\frac{1}{1-\alpha}}$ with $0 < \bar{X}_2 < \bar{X}_1 < 1$. For each fixed $\alpha, \sigma, w$, we temporally assume the existence of a Nash equilibrium of the pure-strategy type, that is, $\mu$ is a Dicar measure. Then, it is concentrated at $x = \bar{X}_2$ if $\bar{X} \leq \bar{X}_2$ and $x = \bar{X}_1$ if $\bar{X} \geq \bar{X}_1$. No optimal pure strategy may exist if $\bar{X} \in (\bar{X}_1, \bar{X}_2)$ and we expect that if an optimal strategy exists, then it should be of a mixed type. The computational experiment implies the existence of Nash equilibria other than the pure-strategy type discussed here and are selected by the classical replicator dynamic. For instance, in the case $\bar{X} \geq \bar{X}_1$, the utility by the mixed strategy almost coincides with that having the optimal pure strategy (see **Figures A1-A3** in the **Appendix**). Such a difficulty is due to the discontinuity of the integrand in the energy cost in (54). The computational resolution used in this subsection is the same as that in the previous one.

The parameter values of utility (54) are set to $\alpha = 0.5$ and $\sigma = w = 1.25$. We first use the homogeneous terminal condition with the terminal time $T = 250$. We examined the classical replicator case (power case with $q = 1$) and exponential case ($q = 1.5$) but present the results for the former case here because they are similar. **Table 7** compares the GRD and MFG for different values of discount rate $\delta$ at time $t = T/2$, again suggesting convergence from the latter to the former although the convergence rate is smaller than the potential cases possibly due to the utility with lower regularity in the present case. Thus, we computationally verify the ansatz that the GRD also arises as a myopic perception limit of the MFG.

The probability density of the strategy distribution is concentrated around a point larger than the threshold $\bar{X} = 0.5$ when no terminal gain is observed, suggesting that the river environment is threatened in this case. To determine the influence of the terminal condition on the probability density, we examine the MFG using the non-homogeneous terminal condition of the form $\Psi(x) = \bar{\Psi}(1-x)$ with a constant $\bar{\Psi} > 0$, with which players select smaller energy use, particularly near the terminal time. We set the terminal time as $T = 100$. The discount rate is $\delta = 0.01$, which is sufficiently small such that the terminal conditions affect the solution profile. Using a larger $\delta$ would be meaningless in this case because the influence of the non-homogenous terminal condition is localized only at the terminal time.

[1] The data available at Water Information System: http://www1.river.go.jp/cgi-bin/SrchDamData.exe?ID=1368041150050&KIND=1&PAGE=0

**Figure 4** compares the probability densities under different terminal conditions, demonstrating that using a larger terminal gain leads to a more concentrated probability density near the boundary $x = 0$. Employing the terminal gain therefore serves as an incentive for the small energy use. This type of policy-making cannot be analyzed using a myopic evolutionary game model where players make decisions based only on the current state. On this basis, the MFG model is more advantageous. Nevertheless, **Figure 4** implies that the actions of the players are not modified for a relatively small time $t$. Strengthening incentives of saving energy through designing utility $U$ in the objective function may resolve this issue. Transformation of society toward sustainable development should then be advanced such that we can find a compromise point between environment and human lives.

**Table 7.** Comparison of the probability density and value function between the generalized classical replicator dynamic and MFG for the energy management case.

| | Probability density | | Value function | |
|---|---|---|---|---|
| $\delta$ | Err | CV | Err | CV |
| $10^{-2}$ | 1.0.E+00 | | 2.8.E-01 | |
| $10^{-1}$ | 1.2.E-01 | 9.37.E-01 | 3.2.E-01 | -6.03.E-02 |
| $10^{0}$ | 1.7.E-02 | 8.43.E-01 | 1.1.E-01 | 4.72.E-01 |
| $10^{1}$ | 2.1.E-03 | 9.06.E-01 | 1.5.E-02 | 8.43.E-01 |
| $10^{2}$ | 4.4.E-04 | 6.77.E-01 | 1.6.E-03 | 9.79.E-01 |

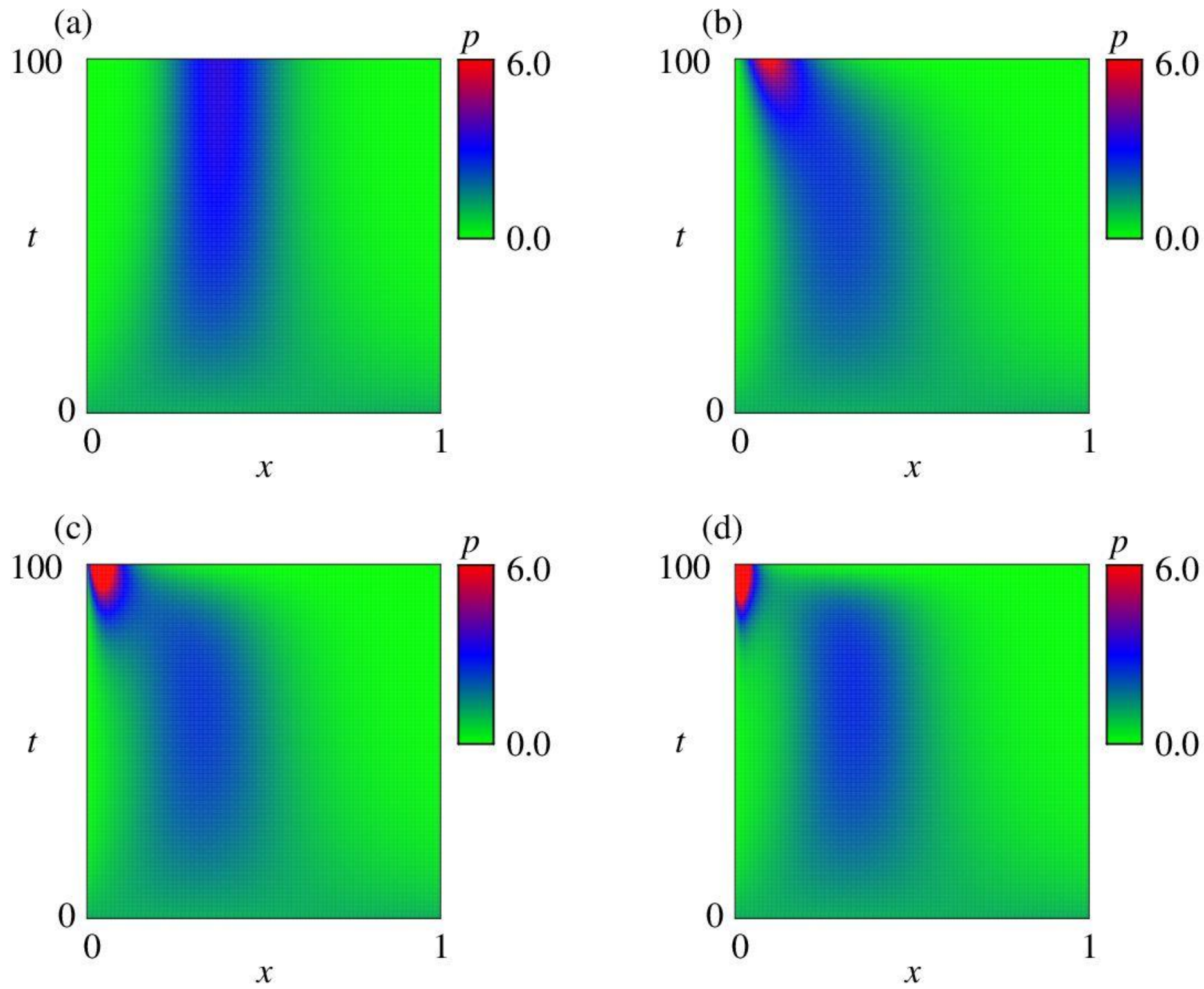

**Figure 4.** Probability densities $p$ under the non-homogeneous terminal condition with different values of $\bar{\Psi}$: (a) $\bar{\Psi} = 0$, (b) $\bar{\Psi} = 1$, (c) $\bar{\Psi} = 2$, and (d) $\bar{\Psi} = 4$.

## 5. Conclusion

An MFG-based explanation of the GRD was proposed, along with computational results. The proposed approach was based on the large-discount limit of an MFG, whose cost function is obtained from an inverse control approach. The transition rate in the GRD was explicitly related to the MFG. We presented a structure-preserving finite difference method for computing evolutionary and MFG models. Numerical solutions to these models were successfully obtained and the convergence between their solutions was computationally analyzed. The computational results suggested that the solutions to the MFG converge to those of the corresponding GRD as the discount rate increases. We consider that the presented computational results are still not directly applicable to real-world problems because they involve a number of simplifications. Nevertheless, the computational examples demonstrated in this paper would help policy-making and understand its consequence depending on foresights of policy-makers. At least, our results suggested that the timescale of myopicity would be critical.

The MFG in this study has several difficulties, such as the non-monotonicity of the operators and the non-separable Hamiltonian. Tackling these issues from a mathematical perspective is necessary to better understand the convergence of the MFG to an evolutionary game. The proposed inverse control approach can be applied to other evolutionary game models not included in this study, such as the higher-order evolutionary game model, as a decision-making model with inertia [76]. An issue to be resolved in this case is the increased number of problem dimensions, whose numerical computation would be much more difficult than that in this study. However, a higher-order model has been applied to machine learning [77] and is therefore worth investigating. Another interesting issue is the investigation of the difference in stationary equilibria between the evolutionary and MFG models [78], such that long-term behaviors of their solutions are better understood.

**Acknowledgements:** This study was supported by the Japan Society for the Promotion of Science: 22K14441 and 22H02456, and Nippon Life Insurance Foundation (Environmental Research Grant for Young Researchers in 2024 No. 24).
**Funding:** This study was supported by the Japan Society for the Promotion of Science: 22K14441 and 22H02456.
**Data availability statement:** The data will be made available upon reasonable request to the corresponding author.
**Competing interests:** The author has no competing interests.
**Declaration of generative AI in scientific writing:** The authors did not use generative AI technology to prepare this manuscript.

## Appendices

### A.1 Proofs

***Proof of Proposition 1***

The proof is a direct consequence of the Proof of Theorem 1 in Cheung [11], along with the boundedness and Lipschitz continuity conditions (2) and (3) which consider a pairwise comparison dynamic in a more abstract setting. Conditions (2) and (3) are required to apply the generalized Picard–Lindelöf theorem (Theorem 3.A in Zeidler [52]) to (12) because this theorem applies to a Banach space with a bounded norm, which is $\mathcal{M}_2(\Omega)$ in our case. The strategy of the Proof of Theorem 1 in Cheung [11] first shows the unique existence of a solution $(\mu_t)_{t\geq 0}$ such that $\mu_t \in \mathcal{M}_2(\Omega)$ at each $t \geq 0$, and then shows $\mu_t \in \mathcal{P}(\Omega)$ based on the conservation of probability encoded in pairwise comparison dynamics.

To use the strategy explained above, we should check the following condition by introducing the transition rate depending on $C$ : if $\mu, \nu \in \mathcal{M}_2(\Omega)$, then a constant $K_1 > 0$ exists such that

$$\sup_{x,y\in\Omega} \rho(y,x,\mu) = \sup_{x,y\in\Omega} C\left(U(x,\mu) - U(y,\mu)\right) \leq K_1 \tag{56}$$

and

$$\sup_{x,y\in\Omega} \left|\rho(x,y,\mu) - \rho(x,y,\nu)\right| \leq K_1 \left\|\mu - \nu\right\|. \tag{57}$$

The first condition (56) is checked as

$$\sup_{x,y\in\Omega} C\left(U(x,\mu) - U(y,\mu)\right) \leq \sup_{x,y\in\Omega} C\left(2\bar{U}\right) = C\left(2\bar{U}\right). \tag{58}$$

The second condition (57) follows from

$$\begin{aligned} \sup_{x,y\in\Omega} \left|\rho(x,y,\mu) - \rho(x,y,\nu)\right| &= \sup_{x,y\in\Omega} \left|C\left(U(x,\mu) - U(y,\mu)\right) - C\left(U(x,\nu) - U(y,\nu)\right)\right| \\ &\leq \sup_{x,y\in\Omega} K_1'\left(\left|U(x,\mu) - U(x,\nu)\right| + \left|U(y,\mu) - U(y,\nu)\right|\right) \\ &\leq 2KK_1'\left\|\mu - \nu\right\| \end{aligned} \tag{59}$$

with a constant $K_1' > 0$ that exists because of the local Lipschitz continuity of the $C$. Setting $K_1 = \max\left\{C\left(2\bar{U}\right), 2KK_1'\right\}$ is sufficient.

□

***Proof of Proposition 2***

The proof is based on the nonlocal degenerate ellipticity and the following direct calculations with an induction argument starting from $i = I$ : for each $i = 0,1,2,...,I-1$ and $j = 1,2,3,...,J$,

$$\begin{aligned} \Phi_{i,j} &= (1-\delta\Delta t)\Phi_{i+1,j} + \Delta t \sum_{k=1}^{J} p_{i,k} \left( \int_0^{(\Phi_{i+1,k} - \Phi_{i+1,j})_+} C(w)\mathrm{d}w \right) \Delta x + \Delta t \delta U_{i,j} \\ &\geq (1-\delta\Delta t)\Phi_{i+1,j} + \Delta t \delta U_{i,j} \\ &\geq (1-\delta\Delta t)\Phi_{i+1,j} - K\delta\Delta t \\ &\geq (1-\delta\Delta t) \min_{m=1,2,3,...,J} \Phi_{i+1,m} - K\delta\Delta t \end{aligned} \tag{60}$$

and

$$\begin{aligned}
\Phi_{i,j} &= \left(1-\delta\Delta t\right)\Phi_{i+1,j} + \Delta t\sum_{k=1}^{J} p_{i,k}\left(\int_0^{\left(\Phi_{i+1,k}-\Phi_{i+1,j}\right)_+} C\left(w\right)\mathrm{d}w\right)\Delta x + \Delta t\delta U_{i,j} \\
&\leq \left(1-\delta\Delta t\right)\Phi_{i+1,j} + \Delta t\sum_{k=1}^{J} p_{i,k}\Delta x L_C\left(\Phi_{i+1,k}-\Phi_{i+1,j}\right)_+ + \delta\Delta t K \\
&\leq \left(1-\delta\Delta t\right)\Phi_{i+1,j} + \Delta t\sum_{k=1}^{J} p_{i,k}\Delta x L_C\left(\max_{m=1,2,3,\ldots,J}\Phi_{i+1,m}-\Phi_{i+1,j}\right)_+ + \delta\Delta t K \\
&= \left(1-\delta\Delta t\right)\Phi_{i+1,j} + \Delta t\sum_{k=1}^{J} p_{i,k}\Delta x L_C\left(\max_{m=1,2,3,\ldots,J}\Phi_{i+1,m}-\Phi_{i+1,j}\right) + \delta\Delta t K \quad . \qquad (61) \\
&= \left(1-\delta\Delta t\right)\Phi_{i+1,j} + L_C\Delta t\left(\max_{m=1,2,3,\ldots,J}\Phi_{i+1,m}-\Phi_{i+1,j}\right) + \delta\Delta t K \\
&= \left(1-\left(\delta+L_C\right)\Delta t\right)\Phi_{i+1,j} + L_C\Delta t\max_{m=1,2,3,\ldots,J}\Phi_{i+1,m} + \delta\Delta t K \\
&\leq \left(1-\left(\delta+L_C\right)\Delta t\right)\max_{m=1,2,3,\ldots,J}\Phi_{i+1,m} + L_C\Delta t\max_{m=1,2,3,\ldots,J}\Phi_{i+1,m} + \delta\Delta t K \\
&= \left(1-\delta\Delta t\right)\max_{m=1,2,3,\ldots,J}\Phi_{i+1,m} + \delta\Delta t K
\end{aligned}$$

From (60), for $i = 0,1,2,\ldots,I-1$ we obtain

$$\begin{aligned}
\min_{m=1,2,3,\ldots,J}\Phi_{i,m} &\geq \left(1-\delta\Delta t\right)^{I-i}\min_{m=1,2,3,\ldots,J}\Phi_{I,m} - K\delta\Delta t\sum_{n=0}^{I-i-1}\left(1-\delta\Delta t\right)^n \\
&\geq \left(1-\delta\Delta t\right)^{I-i}\min_{x\in\Omega}\Psi\left(x\right) - K\delta\Delta t\sum_{n=0}^{I-i-1}\left(1-\delta\Delta t\right)^n \\
&\geq \left(1-\delta\Delta t\right)^{I-i}\min_{x\in\Omega}\Psi\left(x\right) - K\delta\Delta t\frac{1}{1-\left(1-\delta\Delta t\right)} \qquad (62) \\
&= \left(1-\delta\Delta t\right)^{I-i}\min_{x\in\Omega}\Psi\left(x\right) - K \\
&\geq \min_{x\in\Omega}\Psi\left(x\right) - K
\end{aligned}$$

and hence, for $i = 0,1,2,\ldots,I$ ,

$$\min_{m=1,2,3,\ldots,J}\Phi_{i,m} \geq \min_{x\in\Omega}\Psi\left(x\right) - K \ . \qquad (63)$$

Similarly, from (61) we obtain the inequality

$$\max_{m=1,2,3,\ldots,J}\Phi_{i,m} \leq \max_{x\in\Omega}\Psi\left(x\right) + K \quad i = 0,1,2,\ldots,I \ . \qquad (64)$$

The proof completes by setting $K_2 = \max\left\{\left|\max_{x\in\Omega}\Psi\left(x\right)+K\right|,\left|\min_{x\in\Omega}\Psi\left(x\right)-K\right|\right\}$, which is independent of $\Delta t, \Delta x, \delta$ .

□

***Proof of Proposition 3***

The proof of (2.a) is owing to the following direct calculation along with an induction argument:

$$\begin{aligned}
\sum_{j=1}^{J} p_{i,j}\Delta x &= \sum_{j=1}^{J} p_{i-1,j}\Delta x + \Delta x\sum_{j=1}^{J}\Delta x\sum_{k=1}^{J}\left(C\left(\Phi_{i,j}-\Phi_{i,k}\right)-C\left(\Phi_{i,k}-\Phi_{i,j}\right)\right)p_{i-1,k}p_{i-1,j} \\
&= \sum_{j=1}^{J} p_{i-1,j}\Delta x + \left(\Delta x\right)^2\sum_{j=1}^{J}\sum_{k=1}^{J}\left(C\left(\Phi_{i,j}-\Phi_{i,k}\right)p_{i-1,j}p_{i-1,k} - C\left(\Phi_{i,k}-\Phi_{i,j}\right)p_{i-1,k}p_{i-1,j}\right) . \qquad (65) \\
&= \sum_{j=1}^{J} p_{i-1,j}\Delta x
\end{aligned}$$

The proof of (2.b) is owing to (2.a) and the following inequality with an induction argument:

$$
\begin{aligned}
p_{i,j} &= p_{i-1,j}\left(1+\Delta t\Delta x\sum_{k=1}^{J}\left(\left(C\left(\Phi_{i,j}-\Phi_{i,k}\right)-C\left(\Phi_{i,k}-\Phi_{i,j}\right)\right)\right)p_{i-1,k}\right) \\
&\geq p_{i-1,j}\left(1-\Delta t 2C\left(2K_2\right)\sum_{k=1}^{J}p_{i-1,k}\Delta x\right) \\
&\geq p_{i-1,j}\left(1-2C\left(2K_2\right)\Delta t\right)
\end{aligned}
\quad . \qquad (66)
$$

The proof of (2.c) follows from (50) that the right-hand side is proportional to $p_{i-1,j}$.

□

**A.2 Supplementary to Section 4.2**

We first investigate convergence of numerical solutions in $I, J$. We examine the convergence with the parameters $\left(I,J\right)=\left(2500\times 2^{n-1},50\times 2^{n-1}\right)$ ($n=1,2,3,4,5$). The error is measured at time $t=T/2$, which is investigated in the main text. We regard the numerical solution with $n=5$ as the reference solution because the exact solution is not available. For the MFG in the exponential case ($q=1.5$) with concave $f$, **Table A1** summarizes the convergence of numerical solutions to the reference one as $n$ increases. The computational error measured in the maximum norm of the probability density is at most $O\left(10^{-3}\right)$ with $n=3$, which is the resolution employed in the main text. The error is considered sufficiently small.

Next, we investigate the convergence in $\delta$ for the logarithmic and negative exponential cases as summarized in **Tables A2–A5**. The convergence is again at the first order in $\delta^{-1}$. The logarithmic case successfully obtains converged numerical solutions for convex $f$ with a small $\delta$, whereas the other cases fail. This is due to the weaker nonlinearity in this case than the others.

**Table A1.** Numerical convergence of the probability densities of the MFG for the exponential case ($q = 1.5$) with concave $f$. The convergence rate (CV) at the resolution $n$ against $n-1$ is evaluated as $\log_2\left(\mathrm{Err}_{n-1} / \mathrm{Err}_n\right)$ where the subscript represents the resolution.

| | $\delta = 0.01$ | | $\delta = 1$ | |
|---|---|---|---|---|
| $n$ | Err | CV | Err | CV |
| 1 | 6.8.E-04 | | 1.8.E-03 | |
| 2 | 2.8.E-04 | 1.28.E+00 | 7.7.E-04 | 1.25.E+00 |
| 3 | 1.0.E-04 | 1.49.E+00 | 2.0.E-04 | 1.94.E+00 |
| 4 | 4.0.E-05 | 1.32.E+00 | 9.0.E-05 | 1.15.E+00 |

**Table A2.** Same as **Table 1** for the logarithmic case ($q = 0.5$).

| | Probability density | | Value function | |
|---|---|---|---|---|
| $\delta$ | Err | CV | Err | CV |
| $10^{-2}$ | 8.8.E-01 | | 1.8.E-01 | |
| $10^{-1}$ | 9.9.E-02 | 9.48.E-01 | 6.6.E-02 | 4.27.E-01 |
| $10^{0}$ | 1.0.E-02 | 9.78.E-01 | 1.1.E-02 | 7.61.E-01 |
| $10^{1}$ | 1.0.E-03 | 1.00.E+00 | 1.3.E-03 | 9.63.E-01 |
| $10^{2}$ | 1.1.E-04 | 9.76.E-01 | 1.2.E-04 | 1.01.E+00 |

**Table A3.** Same as **Table 1** for the negative exponential case ($q = 2$).

| | Probability density | | Value function | |
|---|---|---|---|---|
| $\delta$ | Err | CV | Err | CV |
| $10^{-2}$ | 1.9.E+00 | | 2.0.E-01 | |
| $10^{-1}$ | 2.1.E-01 | 9.74.E-01 | 1.3.E-01 | 2.01.E-01 |
| $10^{0}$ | 2.1.E-02 | 9.95.E-01 | 3.7.E-02 | 5.32.E-01 |
| $10^{1}$ | 2.1.E-03 | 1.00.E+00 | 4.8.E-03 | 8.83.E-01 |
| $10^{2}$ | 2.1.E-04 | 9.94.E-01 | 5.0.E-04 | 9.86.E-01 |

**Table A4.** Same as **Table A2** for convex $f$.

| | Probability density | | Value function | |
|---|---|---|---|---|
| $\delta$ | Err | CV | Err | CV |
| $10^{-2}$ | 5.8.E+00 | | 2.9.E-01 | |
| $10^{-1}$ | 7.8.E-01 | 8.71.E-01 | 5.5.E-02 | 7.16.E-01 |
| $10^{0}$ | 7.9.E-02 | 9.93.E-01 | 9.7.E-03 | 7.55.E-01 |
| $10^{1}$ | 7.9.E-03 | 1.00.E+00 | 1.0.E-03 | 9.71.E-01 |
| $10^{2}$ | 8.0.E-04 | 9.95.E-01 | 9.5.E-05 | 1.04.E+00 |

**Table A5.** Same as **Table A3** for convex $f$.

| | Probability density | | Value function | |
|---|---|---|---|---|
| $\delta$ | Err | CV | Err | CV |
| $10^{-2}$ | CF | | CF | |
| $10^{-1}$ | CF | | CF | |
| $10^{0}$ | 2.3.E-01 | | 3.6.E-02 | |
| $10^{1}$ | 2.3.E-02 | 1.00.E+00 | 4.7.E-03 | 8.84.E-01 |
| $10^{2}$ | 2.3.E-03 | 1.00.E+00 | 4.9.E-04 | 9.88.E-01 |

### A.3 Supplementary to Section 4.3

Here, we investigate the utility of energy management (54) in the stationary state. For a pure strategy, utility $U$ becomes a function of $x \in \Omega$ and is expressed as

$$U(x) = \frac{1}{\alpha} x^{\alpha} - \sigma x \times \begin{cases} 1 & (x < \bar{X}) \\ 1 + w & (x \geq \bar{X}) \end{cases} \tag{67}$$

with the derivative

$$\frac{\mathrm{d}}{\mathrm{d}x} U(x) = x^{\alpha - 1} - \sigma \times \begin{cases} 1 & (x < \bar{X}) \\ 1 + w & (x \geq \bar{X}) \end{cases}. \tag{68}$$

As shown in **Figure A1**, if $\bar{X} < \bar{X}_2$ then $U$ is locally maximized at $x = \bar{X}_2 = \left(\sigma(1+w)\right)^{-\frac{1}{1-\alpha}}$ with the supremum but not the maximum at $x = \bar{X}$, if $\bar{X} > \bar{X}_1$ then $U$ is globally maximized at $x = \bar{X}_1 = \sigma^{-\frac{1}{1-\alpha}}$, and if $\bar{X}_2 < \bar{X} < \bar{X}_1$, then only a supremum but not the maximum exists at $x = \bar{X}$. For $\bar{X} < \bar{X}_2$, a better pure strategy than that concentrated at $x = \bar{X}_2$ exists.

**Figure A2** shows the computational results of the classical replicator dynamics with $\Delta t = 0.1$ and $\Delta x = 0.01$ at time $t = 4000$. The initial condition is a uniform distribution on $\Omega$. The results suggest the following equilibrium in this case: if $\bar{X} < \bar{X}_2$, then the selected equilibrium is the pure strategy concentrated at $x = \bar{X}_2$, with which the density effect is activated ($\int_{[\bar{X},1]} \mu(\mathrm{d}y) > 0$) and some mixed strategy if $\bar{X}_2 < \bar{X} < \bar{X}_1$. If $\bar{X} > \bar{X}_1$, the utilities at equilibrium are close to each other for the computed and theoretical values, where the latter is assumed to be a pure strategy. This finding suggests the existence of an optimal mixed strategy that maximizes utility. Nevertheless, it should be noted that the convergence of replicator dynamics has been theoretically suggested to be slow and is almost reciprocal of time [79]. **Figure A3** of the average utility $\bar{U} = \int_{\Omega} U(x,\mu)\mu(\mathrm{d}x)$ for different time instances implies that terminating the computation earlier may lead to a biased result, although this may not be critical in this case. This figure also suggests that the optimal strategy switches at $\bar{X} = \bar{X}_1, \bar{X}_2$ as theoretically predicted, and moreover the consistency with **Figure A1** for the cases $\bar{X} > \bar{X}_1$ and $\bar{X} < \bar{X}_2$.

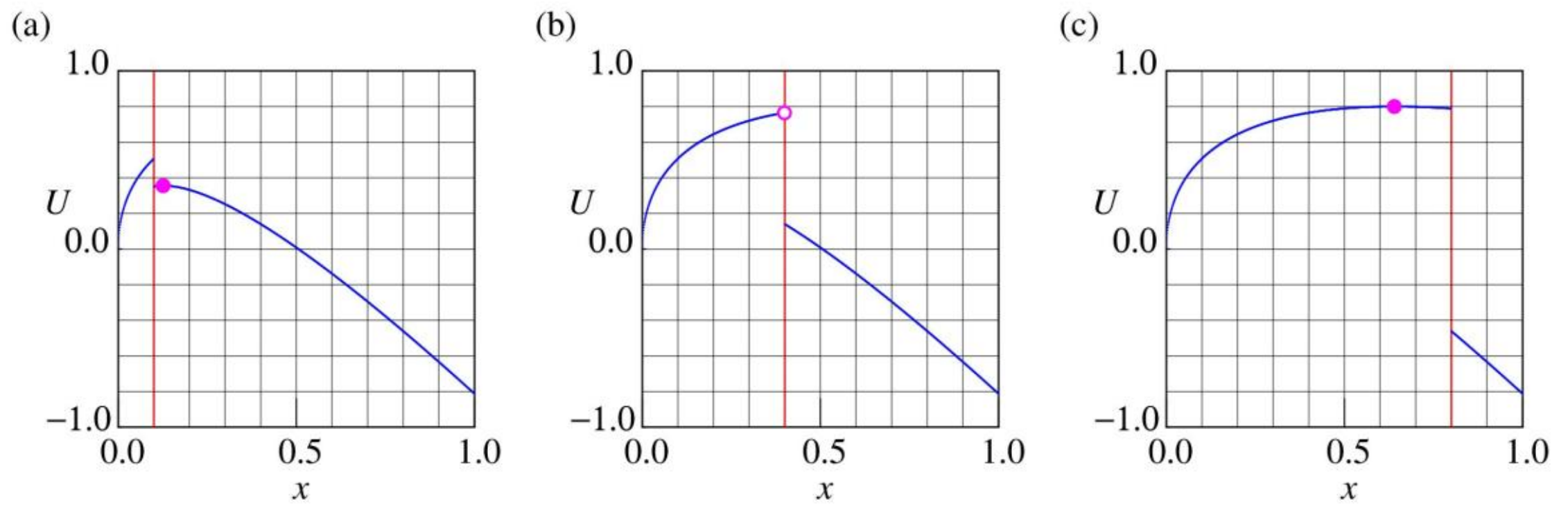

**Figure A1.** Functional shape of $U = U(x)$: (a) $\bar{X} = 0.1 < \bar{X}_2$, (b) $\bar{X}_2 < \bar{X} = 0.4 < \bar{X}_1$, and (c) $\bar{X} = 0.8 > \bar{X}_1$. Here, $\bar{X}_1 = 0.640$ and $\bar{X}_2 = 0.126$. The red line represents discontinuity of $U$, filled magenta circle the local maximum, and unfilled magenta circle supremum but not maximum. The $U$ values of the circles are (a) 0.356, (b) 0.765, (c) 0.800.

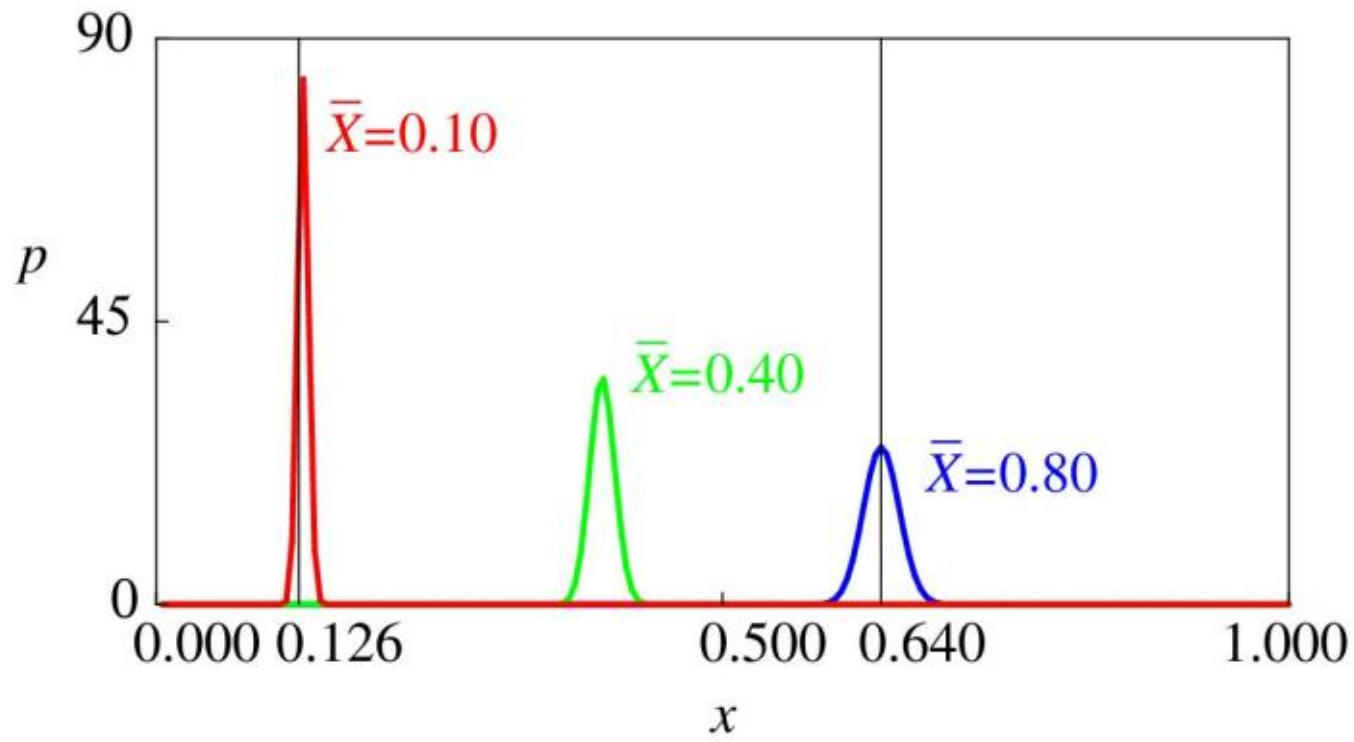

**Figure A2.** Solutions to the classical replicator dynamic for different values of at time $t = 4000$.

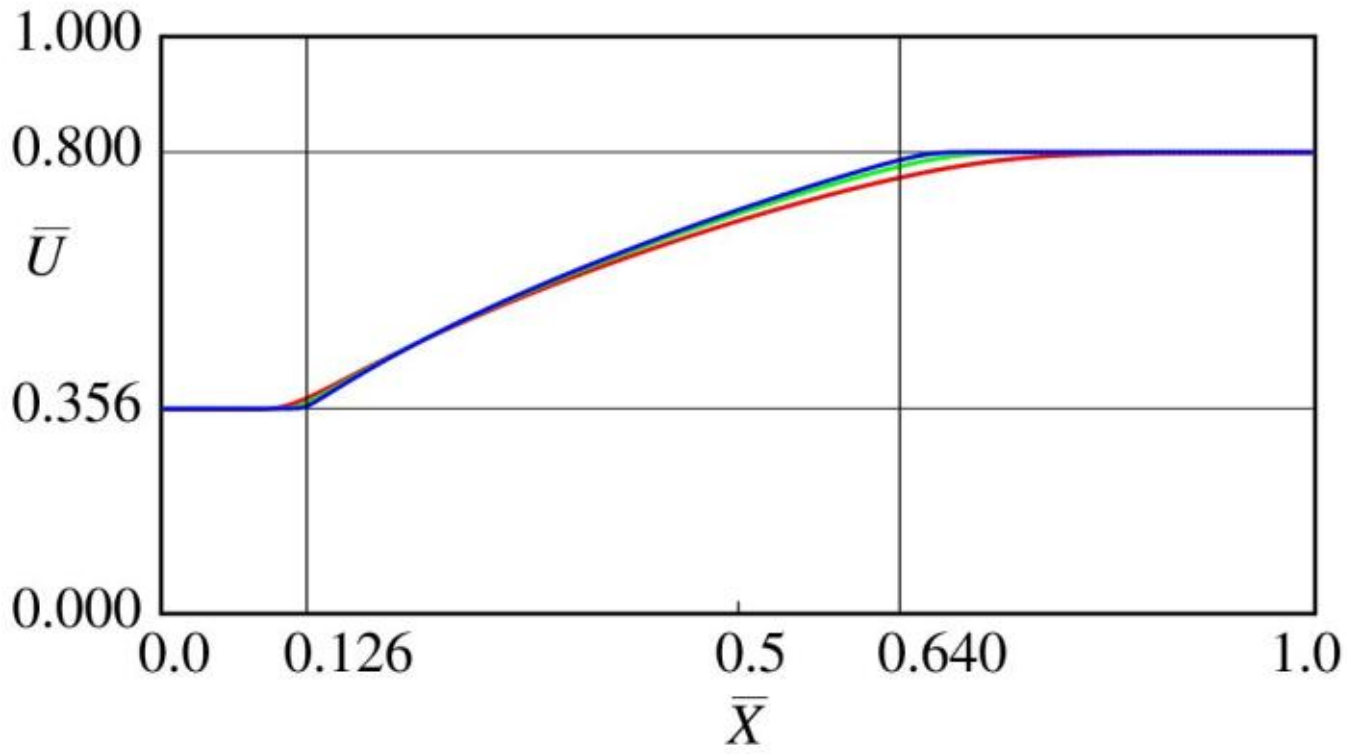

**Figure A3.** Computed utilities $\bar{U}$ at different time instances: $t = 250$ (red), $t = 1000$ (green), and $t = 4000$ (blue). The two horizontal lines correspond to the global maximum for $\bar{X} < \bar{X}_2 = 0.126$ ($\bar{U} = 0.356$) and the local maximum for $\bar{X} > \bar{X}_1$ ($\bar{U} = 0.800$).